\documentclass{article}
\usepackage[T1]{fontenc}
\usepackage{amsmath, amssymb, amsfonts, amsthm}
\usepackage{geometry}
\usepackage{bm}
\usepackage{graphicx}
\usepackage{textcomp}
\usepackage{mathrsfs}
\usepackage{url}
\usepackage{mathtools}
\allowdisplaybreaks[4]
\usepackage{tikz-cd}
\usepackage{hyperref}
\usepackage{enumitem}
\usepackage{bm}
\usepackage[all]{xy}

\newtheorem*{TECNC}{The  coarse Novikov Conjecture}
\newtheorem*{TEBCC}{The  coarse Baum-Connes Conjecture}
\theoremstyle{plain}
\newtheorem{theorem}{Theorem}[section]
\newtheorem{lemma}[theorem]{Lemma}

\newtheorem{proposition}[theorem]{Proposition}

\newtheorem{corollary}[theorem]{Corollary}
\theoremstyle{definition}

\newtheorem{definition}[theorem]{Definition}

\renewcommand{\epsilon}{\varepsilon}
\theoremstyle{remark}

\newtheorem{remark}[theorem]{Remark}

\DeclareMathOperator{\id}{id}

\def\A{\mathcal A}

\def\B{\mathcal B}
\def\K{\mathcal K}
\def\H{\mathcal H}

\def\IR{\mathbb R}
\def\IN{\mathbb N}
\def\IC{\mathbb C}
\def\Cl{\textup{Cliff}_{\IC}}
\def\S{\mathcal S}

\begin{document}
	\title{$\ell^p$-Busemann space and the coarse Novikov conjecture}
	\vskip 16mm
	\author{Chao Tong and\ Xianjin Wang\footnote{The authors are supported in part by NSFC (No. 12071183 and No. 11771061).}
	}

	\date{\today}                                  
	\maketitle
	\begin{abstract}
	We introduce the notion of an \(\ell^p\)-Busemann space, which is a complete Busemann convex geodesic metric space whose tangent cones admit uniformly controlled homogeneous bi-Lipschitz embeddings into \(\ell^p\)-spaces. An \(\ell^p\)-Busemann space is called admissible if it contains dense sequences of closed convex subsets isometric to certain Finsler manifolds. We then prove that the rational coarse Novikov conjecture holds for proper metric spaces with bounded geometry that admit coarse embeddings into admissible \(\ell^p\)-Busemann spaces.
	\end{abstract}
	
	\par\noindent {\bfseries Keywords.} Metric space, Busemann convex space, Higher index theory, Coarse Novikov conjecture.

	\section{Introduction}
	The Novikov conjecture is a fundamental problem in the topology of manifolds. It asserts the homotopy invariance of certain higher signatures of a manifold $M$. By methods of noncommutative geometry, the Novikov conjecture is a consequence of the injectivity of the Baum-Connes assembly map. Let $G=\pi_1(M)$ be a countable discrete group, and let $\underline EG$ denote
	the universal proper $G$-space. The Baum--Connes assembly map is a
	homomorphism
	\[
	\mu_G:
	\varinjlim_{Z\subseteq \underline EG,\; G\text{-compact}}
	KK_*^G(C_0(Z),\mathbb C)
	\to
	K_*(C_r^*(G)),
	\]
	 The Baum--Connes
	conjecture asserts that $\mu_G$ is an isomorphism. Its injectivity
	part is usually referred to as the strong Novikov conjecture. The rational strong Novikov conjecture states that the assembly map is injective when tensoring with $\mathbb{Q}$. For general metric spaces, there is also an coarse Baum-Connes conjecture, which states that the assembly map $$\mu :KX_*(X)\to K_*(C^*(X))$$ 
	is an isomorphism. The left-hand side is the coarse K-homology groups of a metric space, and the right-hand side is the  Roe algebra of the metric space. The Baum-Connes conjecture for $G$ is the coarse Baum-Connes conjecture for $X=G$ as a metric space. The rational strong coarse Novikov conjecture implies the Novikov conjecture and the Gromov-Lawson-Rosenberg conjecture regarding the existence of metrics of positive scalar curvature on aspherical closed manifolds.

	There are several models for K-homology. One of them is the localization algebra approach introduced by Guoliang Yu. \cite{yu1997localization}. It has the advantage that the assembly map is given by an explicit $*$-homomorphism, the evaluation map. Using this approach, Yu proved the Novikov conjecture for groups with finite asymptotic dimension \cite{yu1998novikov}, and also the coarse Baum-Connes conjecture for spaces admitting a coarse embedding into Hilbert spaces \cite{yu2000coarse}. However, coarse embeddability into Hilbert space is not a necessary condition for the injectivity of the coarse assembly map. This becomes particularly significant for expander graphs, which cannot coarsely embed into a Hilbert space, and certain expanders give counterexamples to the coarse Baum-Connes conjecture \cite{9b71e94cd22848fc8803555c49b2b1e9}. Since then, various developments have been made for the coarse Novikov conjecture for spaces which coarsely embed into suitable spaces. G. Kasparov and G. Yu proved the coarse Novikov conjecture for spaces uniformly coarsely embedding into uniformly convex Banach spaces \cite{KASPAROV20061}, and they also proved the Novikov conjecture for spaces coarsely embedding into Banach spaces with property $(H)$ \cite{kasparov2012novikov}. Jianchao Wu and others considered spaces coarsely embedding into Hilbert-Hadamard spaces \cite{gong2023novikov},\cite{gong2021novikov}, and concluded that certain diffeomorphism groups of manifolds satisfy the Novikov conjecture. A major development beyond coarse embeddability is the notion of fibred coarse embedding into a Hilbert space, introduced by Chen, Wang, and Yu \cite{CHEN201388}. This condition guarantees that the maximal coarse Baum–Connes conjecture holds. Subsequently, Finn-Sell used the coarse groupoid approach to establish injectivity of the coarse assembly map under fibred coarse embeddability assumptions, even in situations where ordinary coarse embeddability into Hilbert space fails \cite{FINNSELL20143758}. 
	
	There has been extensive research on the Novikov conjecture, here we do not list all of it. Our paper is inspired by Wu's Hilbert-Hadamard spaces \cite{gong2023novikov},\cite{gong2021novikov}, Yu's Property $(H)$ \cite{kasparov2012novikov}, and Wang's $\ell^p$-spaces approach \cite{guo2024bott}. In Wu's work, a Hilbert-Hadamard space is a CAT(0) geodesic metric space, and every point has a tangent cone which isometrically embeds into a Hilbert space, allowing Clifford operators and Bott maps to be defined. We would like to generalize their methods to metric spaces such that every point has a tangent cone which  embeds into an $\ell^p$ space, so that we can apply Wang's method, using the Mazur map (a homeomorphism from the unit ball of $\ell^p$ to $\ell^2$), to define the Clifford operator and Bott map. However, $\ell^p$ is not CAT$(0)$ unless $p=2$, so we need a more general notion of metric spaces with nonpositive curvature, namely the \emph{Busemann convex} spaces (or spaces with \emph{Busemann nonpositive curvature}). Thus, we introduce the \emph{$\ell^p$-Busemann space}: a geodesic metric space with nonpositive curvature in the Busemann sense, such that tangent cones admit bi-Lipschitz embeddings into an $\ell^p$ space. The change from an isometric embedding in \cite{gong2021novikov} to a bi-Lipschitz embedding is essential for our construction. We would like to have a Finsler manifold such that the tangent space is a finite-dimensional $\ell^p$ space, however, the standard $\ell^p$ norm does not define a smooth Finsler structure on the manifold. By modifying it to be smooth, the tangent space no longer isometrically embeds, but instead bi-Lipschitz embeds into an $\ell^p$ space. With the help of Mazur map, we can define the Clifford operator and a $C^*$-algebra associated with the $\ell^p$-Busemann space. Under certain additional admissibility conditions, we can construct a Bott map which induces injection on K-theory, and prove the rational coarse Novikov conjecture holds for spaces that coarsely embed into an admissible $\ell^p$-Busemann space.
	Our main theorem states that:
	\begin{theorem}\label{mainthm}
		Let \(X\) be a proper metric space with bounded geometry that admits a coarse embedding into an admissible \(\ell^p\)-Busemann space \(\mathcal M\). Then the rational coarse Novikov conjecture holds for \(X\).
	\end{theorem}
	The paper is organized as follows. In Section 2, we introduce $\ell^p$-Busemann spaces and the coarse Novikov conjecture. In Section 3, we define the twisted Roe algebra and localization algebra, and construct the Bott map from the Roe algebra/localization algebra to the twisted Roe algebra/localization algebra. In Section 4, we prove the twisted assembly map is an isomorphism by a cutting-and-pasting argument, and then prove our main theorem.
	
	\section{Preliminaries}
	\subsection{$\ell^p$-Busemann space}
	
	Let $(M,d)$ be a metric space. A continuous map from the interval $[0,T]$ to $M$ is called a \emph{path}. The \emph{length} of a path $\gamma$, denoted by $\|\gamma\|$, is defined to be
	$$\|\gamma\|=\sup\sum_{i=1}^nd(\gamma(t_{i-1}),\gamma(t_i))$$
	where the supremum is taken over the set of all partitions $0=t_0< \cdots < t_n = T$ of the interval $[0,T]$. We denote $\gamma(t)$ by $\gamma_t$. A path $\gamma:[0,T]\to M$ is called a \emph{geodesic} if $d(\gamma_s,\gamma_t)=\frac{|s-t|}{T}d(\gamma_0,\gamma_T)$ for every $s,t\in [0,T]$. A \emph{geodesic metric space} is a metric space such that any two points can be connected by a geodesic path. A geodesic metric space is called \emph{geodesically complete} if every geodesic segment can be extended indefinitely to a line $\mathbb{R} \to X$. For a general metric space, there are many kinds of notions of curvature that generalize the curvature of a Riemannian manifold. Alexandrov's comparison curvature is the most widely known one. A metric space $X$ is said to be CAT$(k)$ if any triangle in $X$ is thinner than its comparison triangle in the model space $M_k$, the simply connected, complete Riemannian manifold with constant sectional curvature $k$. CAT($k$) spaces with $k\leq 0$ have been well studied, one may refer to \cite{bridson2013metric} for a comprehensive introduction. However, the CAT($0$) condition is too rigid. $\ell^p$ spaces do not satisfy CAT($0$) unless $p=2$. There is a more general notion of \emph{Busemann convexity} or \emph{Busemann nonpositive curvature} for metric spaces, which characterizes curvature by the property that geodesics move further and further away from each other.
	\begin{definition}
		Let $X$ be a geodesic metric space. $X$ is said to have \emph{Busemann nonpositive curvature} or be \emph{Busemann convex} if for any two geodesics $\gamma_1: [0, 1] \to X$ and $\gamma_2: [0, 1] \to X$, the distance function: $$t \longmapsto d(\gamma_1(t), \gamma_2(t))$$ is a convex function for $t \in [0, 1]$. That is, for all $t_1, t_2 \in [0, 1]$ and $\lambda \in [0, 1]$:
		$$d\left(\gamma_1\left((1-\lambda)t_1 + \lambda t_2\right), \gamma_2\left((1-\lambda)t_1 + \lambda t_2\right)\right) \le (1-\lambda)d(\gamma_1(t_1), \gamma_2(t_1)) + \lambda d(\gamma_1(t_2), \gamma_2(t_2))$$
	\end{definition}
	Typical examples include trees and CAT($k$) spaces for $k\leq 0$. An $\ell^p$ space is Busemann convex due to the convexity of the $\ell^p$ norm for \(1<p<\infty\). Another important class of examples is provided by suitable \emph{Berwald-Finsler} manifolds.
	\begin{definition}
		A \emph{Finsler manifold} is a pair $(M, F)$, where $M$ is a smooth manifold and $F: TM \to [0, \infty)$ is a Finsler function (or Finsler norm) satisfying the following conditions:
		\begin{enumerate}
			\item Smoothness: $F$ is $C^\infty$ on the slit tangent bundle $TM_0 = TM \setminus \{0\}$ and continuous on the zero section.
			\item Positive Homogeneity: For any point $x \in M$, any tangent vector $y \in T_x M$, and any scalar $\lambda > 0$: $$F(x, \lambda y) = \lambda F(x, y).$$
			\item Strong Convexity: The fiberwise Hessian of \(F^2\) with respect to the tangent space coordinates is positive definite on the slit tangent bundle $TM_0$, i.e., 
			$$g_{ij}(x, y) = \frac{1}{2} \frac{\partial^2 (F^2)}{\partial y^i \partial y^j}(x, y)$$ is positive-definite for all $(x, y) \in TM_0$.
		\end{enumerate}
	\end{definition}
	
	\begin{definition}
		A Finsler manifold is called \emph{Berwald} if there exists a linear connection on the tangent bundle such that parallel translation preserves the Finsler function $F$.
	\end{definition}
	It is proved in \cite{KRISTALY20061257} that a Berwald Finsler manifold with nonpositive flag curvature is Busemann convex. Let $(X,d)$ be a Busemann convex space and let $x\in X$. Let $\mathcal G_x$ be the set of geodesic germs starting at $x$. More concretely, a representative is a geodesic $\gamma:[0,\varepsilon_\gamma]\to X$ with $\gamma(0)=x$. Two such representatives $\gamma,\eta$ define the same germ, written $\gamma\sim_0\eta$ if they agree on some interval $[0,\varepsilon]$ with $\varepsilon>0$. Form the tangent cone by 
	$$
	\mathcal C_x=\bigl(\mathcal G_x\times[0,\infty)\bigr)/\sim_0,
	$$
	where all pairs $(\gamma,0)$ are identified as a single vertex, denoted by $0$. We write its elements as $[\gamma,s]$. For $[\gamma,s],[\eta,t]\in\mathcal C_x$, define
	$$
	\delta_x([\gamma,s],[\eta,t])=\lim_{r\downarrow0}\frac{d(\gamma(rs),\eta(rt))}{r}.
	$$
	The above limit defines a pseudometric on the pre-tangent cone \cite{bridson2013metric},\cite{GigliNobili2021}. And the metric tangent cone is obtained by identifying points at zero pseudodistance, define
	$$
	u\sim v
	\quad\Longleftrightarrow\quad
	\delta_x(u,v)=0.
	$$
	The metric tangent cone is then the quotient
	$$T_xX=\mathcal C_x/{\sim},\qquad d_{T_xX}([\gamma,s],[\eta,t])=\delta_x([\gamma,s],[\eta,t]).$$
	If $X$ is a CAT$(0)$ space, it is straightforward to check that the metric on the tangent cone coincides with the metric defined by the Euclidean cosine law.
	
	The logarithm function $\log_x: X\to T_xX$ is defined by 
	$$\log_x(y)=[\gamma_{x,y}, d(x,y)],$$
	where $\gamma_{x,y}$ is the unit-speed geodesic from $x$ to $y$. By Busemann convexity, for $0<r\leq 1$,
	$$d(\gamma_{xy}(r d(x,y)),\gamma_{xz}(r d(x,z)))\leq r\,d(y,z).$$
	Dividing by $r$ and letting $r\to 0^+$, we have 
	$$ d_{T_xX}(\log_x(y),\log_x(z))\leq d_X(y,z).$$
	So $\log_x$ is $1$-Lipschitz. Note that if $y,z$ are on the same geodesic $\gamma$ emanating from $x$, then
	$$ d_{T_xX}(\log_x(y),\log_x(z))= d_X(y,z).$$
	A metric space $X$ is said to have a \emph{$C$-bi-Lipschitz embedding} into another metric space $Y$ if there is a constant $C>0$ and a map $f:X \to Y$ such that $$C^{-1}d_X(x,y)\leq d_Y(f(x),f(y))\leq Cd_X(x,y)$$
	\begin{definition}
		An \emph{$\ell^p$-Busemann space} $X$ is a complete,
		Busemann convex geodesic metric space such that
		there exists a constant $C\geq 1$ and, for every
		$x\in X$, a homogeneous $C$-bi-Lipschitz embedding $J_x:T_xX\longrightarrow\ell^p(\mathbb N,\mathbb R)$
		satisfying:
		\begin{enumerate}
			\item $J_x(0)=0$,
			\item $J_x(\lambda v)=\lambda J_x(v),\  \lambda\geq0.$
		\end{enumerate}
	\end{definition}
	
	\begin{definition}
	An $\ell^p$-\emph{Finsler manifold} is a parallelizable Finsler manifold $(M, F)$ such that at each point $x\in M$, the Finsler norm on $TM_x$ is defined by
	$$F_{\epsilon}(x,y)=\left(\sum_{i=1}^n(y_i^2+\varepsilon^2\|y\|_2^2)^{p/2}\right)^{1/p}.$$
	Where $y=\sum_{i=1}^n y_iX_i(x)$, and $X_1,\ldots,X_n$ is the global frame.
\end{definition}
	\begin{remark}
		The function $F(x,y)=\|y\|_p$ does not define a Finsler structure in the strict sense, since $F$ may fail to be $C^\infty$ smooth when $y$ has zero coordinates. However, the function $F_{\epsilon}$ does define a smooth Finsler structure and converges to the $\ell^p$ norm as $\epsilon\to 0$. Let $(V,F_{\epsilon})$ be a finite-dimensional vector space with the norm $F_{\epsilon}$, there may not exist an isometric embedding of $(V,F_{\epsilon})$ into $\ell^p$ spaces, but the identity map $\id:(V,F_{\epsilon}) \to (V,\ell^p)$ is bi-Lipschitz. 
	\end{remark}
	Let $N_d=\{1,2,\dots d\}\subset \IN$. After choosing a suitable linear frame of \(T_xM\), we identify the tangent space with \(\mathbb R^n\) and use the coordinate identity map as a bi-Lipschitz embedding $J_x: T_xM\to \ell^p(N_d,\IR)$:
	$$J_x(y)=y$$ we see that for a $d$-dimensional $\ell^p$-Finsler manifold, the tangent space at each point has a bi-Lipschitz embedding into $\ell^p(N_d,\IR)$, and the embedding is also homogeneous. For notational convenience, we will not distinguish a point $\log_x(y) \in T_xX$ from its image $J_x(\log_x(y))$ in the embedded $\ell^p$ space if it causes no confusion.
	\begin{definition}
		An \emph{$\ell^p$-Finsler-Busemann manifold} is a complete, simply connected, Berwald $\ell^p$-Finsler manifold $(M, F)$ with nonpositive flag curvature. 
	\end{definition}
	Thus, an $\ell^p$-Finsler-Busemann manifold is an $\ell^p$-Busemann space in our sense. For such a manifold, the $\log$ map is the inverse of the ordinary exponential map. The Cartan-Hadamard theorem can be generalized to a Finsler-Cartan-Hadamard theorem, one may refer to \cite{BaoChernShen2000} for a  proof.
	\begin{theorem}
		If $(M,F)$ is a complete, simply connected Finsler manifold of nonpositive flag curvature, then for every $p\in M$,
		$$
		\exp_p:T_pM\longrightarrow M
		$$
		is a $C^1$-diffeomorphism.
	\end{theorem}
	\begin{definition}
		A separable $\ell^p$-Busemann space $M$ is said to be \emph{admissible} if there exists an increasing sequence of subspaces $\{M_n\}_{n\in\mathbb{N}}$ of $M$ such that
		\begin{itemize}
			\item[(1)] each $M_n$ is a convex closed subspace of $M$ and isometric to a finite-dimensional $\ell^p$-Finsler-Busemann manifold,
			\item[(2)] $\bigcup_{n\in\mathbb{N}}M_n$ is dense in $M$.
		\end{itemize}
	\end{definition}
	\subsection{The coarse Novikov conjecture}
	
	In this section we recall some basic facts of the coarse Baum-Connes conjecture. The Roe algebra was introduced by John Roe \cite{roe1993coarse}, it is a $C^*$-algebra that characterizes the large-scale geometry of metric spaces. The K-theory group of the Roe algebra is the receptacle for studying higher index theory of abstract elliptic operators. The coarse Baum-Connes conjecture states that a certain assembly map from the coarse K-homology of a metric space to the K-theory of the Roe algebra is an isomorphism. We will follow the localization algebra approach introduced in \cite{yu1997localization} for a model of $K$-homology. Firstly, we recall some terminology in coarse geometry. 
	
	In this paper we only consider \emph{proper} metric spaces, i.e., metric spaces in which every closed and bounded subset is compact. Let $X, Y$ be proper metric spaces, and $f: X\to Y$ be a (not necessarily continuous) map. Then $f$ is called \emph{proper} if the preimage of any compact set is compact. The map $f$ is said to be a coarse embedding from $X$ to $Y$ if there exist functions $\rho_-,\rho_+:[0,\infty)\to[0,\infty)$ tending to infinity as $t\to \infty$, such that:
	$$\rho_-(d(x_1,x_2))\leq d(f(x_1),f(x_2))\leq \rho_+(d(x_1,x_2)).$$
	For any $r>0$, a subset $X'\subseteq X$ is said to be an \emph{$r$-net} if for any $x\in X$, there exists $y\in X'$ such that $d(y,x)\leq r$.
	Moreover, $f$ is said to be a \emph{coarse equivalence} if $f$ is a coarse embedding and there exists $R>0$ such that the image $f(X)$ is an $R$-net in $Y$.
	
	A discrete metric space $X$ is said to have \emph{bounded geometry} if for any $R>0$, there exists $N>0$ such that $\sup_{x\in X}\#B(x, R)\leq N$. A non-discrete metric space $Y$ is said to have \emph{bounded geometry} if there exists an $r$-net $Y'$ of $Y$ such that $Y'$ has bounded geometry. Since every proper metric space with bounded geometry is coarsely equivalent to a uniformly discrete net of bounded geometry, we may work with such a net whenever convenient.
	\begin{definition}(cf.\cite{willett2020higher},\cite{kasparov2012novikov})
		Let $C_0(X)$ be the $C^*$-algebra of all complex-valued continuous functions on $X$ vanishing at infinity. A Hilbert space $H_X$ is called an \emph{$X$-module} if there is a representation $\varphi:C_0(X)\to\B(H_X)$. $H_X$ is called \emph{admissible} or \emph{ample} if in addition:
		\begin{enumerate}
			\item $\varphi_X(f)\notin \K(\H_X)$ for all nonzero $f\in C_0(X)$,
			\item $\varphi_X$ is non-degenerate in the sense that $\{\varphi_X(f) \xi: f\in C_0(X), \xi\in \H_X\}$ is dense in $\H_X$. 
		\end{enumerate}
	\end{definition}
	
	\begin{definition}
		Let $H_X$ be an $X$-module, and let $T\in\B(H_X)$ be an operator on $H_X$.
		\begin{itemize}
			\item[(1)] The {\em support} of $T$ is defined to be the complement in $X\times X$ of the set of all points $(x, y)\in X\times X$ for which there exist $f, g\in C_0(X)$ such that $fTg=0$ but $f(x)\not= 0, g(y)\not=0$.
			\item[(2)] The {\em propagation} of $T$ is defined to be
			$$\textup{Prop}(T)=\sup\big\{ d(x,y) \mid (x, y)\in \textup{supp}(T) \big\}.$$
			\item[(3)] The operator $T: H_X\to H_X$ is said to be {\em locally compact} if $fT$ and $Tf$ belong to $\K(H_X)$ for all $f\in C_0(X)$, where $\K(H_X)$ is the algebra of all compact operators on $H_X$.
	\end{itemize}\end{definition}
	
	In fact, there is a standard way to construct an admissible $X$-module (see \cite{willett2020higher}).
	Let $N\subset X$ be a countable dense subset of $X$, and let $\H_0$ be an infinite-dimensional separable Hilbert space. Then the Hilbert space
	$$H_X:=\ell^2(N)\otimes \H_0$$ is an admissible $X$-module.

	\begin{definition}
		Let $H_X$ be an admissible $X$-module.
		\begin{itemize}
			\item[(1)] The {\em algebraic Roe algebra}, denoted by $\mathbb C[X]$, is defined to be the $*$-algebra of all locally compact, finite-propagation operators in $\B(H_X)$.
			\item[(2)] The {\em Roe algebra}, denoted by $C^*(X)$, is defined to be the norm completion of $\mathbb{C}[X]$ in $\B(H_X)$.
		\end{itemize}
	\end{definition}
	
	For discrete metric spaces with bounded geometry, the following "X-by-X" matrix characterization of the Roe algebra is more useful:
	
	\begin{definition}[c.f.\cite{yu2000coarse}]
		Let $N\subset X$ be a countable dense subset of $X$, and let $\H_0$ be an infinite-dimensional separable Hilbert space. Let $\mathbb C_{\rm f}[X]$ be the set of all bounded functions $T$ on $N\times N$ such that:
		\begin{itemize}
			\item[(1)] $T(x,y)\in \K(\H_0)$ for all $x,y\in N$,
			\item[(2)] there exists $r>0$ such that $T(x,y)=0$ for all $x,y\in N$ with $d(x,y)>r$,
			\item[(3)] for each bounded set $B\subseteq X$, the set
			$$\{(x,y)\in (N\times N)\cap (B\times B)\mid T(x,y)\ne0\}$$
			is finite,
			\item[(4)] there exists $L>0$ such that
			$$\#\{y\in N\mid T(x,y)\ne0\}\leq L,\quad\#\{y\in N\mid T(y,x)\ne0\}\leq L$$
			for all $x\in N$.
		\end{itemize}
		The algebraic structure on $\mathbb C_{\rm f}[X]$ is given by
		$$(TS)(x,y)=\sum_{z\in N}T(x,z)S(z,y)$$
		$$(T^*)(x,y)=(T(y,x))^*.$$
		$T$ acts on $\delta_x\otimes v\in\ell^{2}(N)\otimes\H_0$ by standard matrix multiplication:
		$$T(\delta_x\otimes v)=\sum_{y\in N}\delta_y\otimes T(y,x)(v),$$
	\end{definition}
	$\mathbb C_{\rm f}[X]$ is a dense subalgebra of $C^*(X)$, and we will use $\mathbb C_{\rm f}[X]$ to replace $\mathbb C[X]$ to define the Roe algebra of $X$. If there exists a coarse equivalence $f:X\to Y$, then there exists a unitary operator $V\in\B(H_X,H_Y)$ covering $f$ such that the map $\operatorname{ad}_V: T\mapsto VTV^*$ gives rise to an isomorphism 
	$$\operatorname{ad}_V:C^*(X)\xrightarrow{\cong} C^*(Y).$$
	One may refer to \cite[Chapter 4]{yu2000coarse} for a detailed construction and proof. As a result, since the identity map $\id: X\to X$ is a coarse equivalence, the Roe algebra $C^*(X)$ does not depend on the choice of the $X$-module. 
	
	Let $X$ be a proper metric space. The $K$-homology group $K_*(X)$ is the dual of the $K$-theory of $X$. It is defined to be $KK_*(C_0(X),\mathbb C)$ by using Kasparov's $KK$-theory \cite{kasparov1988equivariant}. 
	The localization algebra approach to $K_*(X)$ was introduced by G. Kasparov and G. Yu in \cite{kasparov2012novikov}, which is more computable and does not require advanced technical skills in $KK$-theory. One may also refer to \cite{yu1997localization, willett2020higher, qiao2010localization} for an introduction to localization algebras.
	\begin{definition}\label{localization algebra}
		Let $H_X$ be an admissible $X$-module.
		\begin{itemize}
			\item[(1)] The {\em algebraic localization algebra}, denoted by $\mathbb{C}_L[X]$, is defined to be the $*$-algebra of all bounded and uniformly continuous functions $f:[0,\infty)\to\mathbb C_{\rm f}[X]$ such that
			$$\textup{Prop}(f(t))\to0$$
			as $t\to\infty$.
			\item[(2)] The {\em localization algebra}, denoted by $C^*_L(X)$, is defined to be the completion of $\mathbb{C}_L[X]$ under the norm
			$$\|f\|=\sup_{t\in[0,\infty)}\|f(t)\|.$$
	\end{itemize}\end{definition}
	
	$C^*_L(X)$ does not depend on the choice of the admissible $X$-module, one may refer to \cite{kasparov2012novikov} or \cite[Chapter 6]{yu2000coarse} for a proof. The localization approach and the $KK$-theory approach to $K$-homology are connected by the following theorem:
	
	\begin{theorem}[\cite{kasparov2012novikov}]
		If $X$ is a finite-dimensional simplicial polyhedron, then there is a local index map
		$$\textup{Ind}_L:KK_*(C_0(X),\mathbb C)\to K_*(C^*_L(X))$$
		which is an isomorphism.
	\end{theorem}
	
	Next we recall the Rips complexes associated with $X$, which serve as a model of classifying spaces. For each $n\in\mathbb N$, let $S(\mathbb{R}^n)$ be the sphere of radius one in an $n$-dimensional Euclidean space. The metric on $S(\mathbb{R}^n)$ is defined by
	$$d_{S(\mathbb R^n)}(x,y)=\arccos(\langle x,y\rangle),$$
	where $\langle x,y\rangle$ denotes the standard inner product of $\mathbb R^n$. This is the usual Riemannian metric on $S(\mathbb{R}^n)$.
	
	\begin{definition}[Rips complex]
		For each $d\geq 0$, the \emph{Rips complex} of $X$ at scale $d$, denoted $P_d(X)$, consists as a set of all formal sums
		$$z=\sum_{x\in X}t_xx$$
		such that each $t_x$ is in $[0,1]$, $\sum_{x\in X}t_x=1$, and only finitely many $t_x$ are non-zero. The support of $z$ defined by $\textup{supp}(z):=\{x\in X\,|\,t_x\ne0\}$ has diameter at most $d$.
	\end{definition}
	
	Let $F\subseteq X$ be a finite set. Let $\mathbb{R}^F=\{\sum_{z\in F}t_zz\mid t_z\in\mathbb R\text{ and }z\in F\}$ be the finite-dimensional real Euclidean space spanned by $F$. Denote by $X(F)$ the simplex generated by $F$, i.e.,
	$$X(F)=\left\{\sum_{z\in F}t_zz\,\Big|\,t_z\geq 0\text{ and }\sum_{z\in F}t_z=1\right\}\subseteq \mathbb R^F.$$
	Define 
	$$f:X(F)\to\mathcal{S}(\mathbb{R}^F)\quad\text{by}\quad \sum_{z\in F}t_zz\mapsto\sum_{z\in F}\frac{t_z}{\sqrt{\sum_{z\in F}t^2_z}} z$$
	The spherical metric on $X(F)$ is the metric defined by
	$$d_S(x,y)=\frac{2}{\pi}d_{S(\mathbb R^F)}(f(x),f(y)),$$
	and $X(F)$ equipped with this metric is called the \emph{spherical simplex} on $F$. The factor $2/\pi$ is chosen so that $X(F)$ has diameter one.
	
	For any $x,y\in P_d(X)$, a semi-simplicial path $\delta$ between $x$ and $y$ is a finite sequence of points
	$$x=x_0,y_0,x_1,y_1,\cdots,x_n,y_n=y.$$
	The length of $\delta$ is defined to be
	$$l(\delta)=\sum_{i=0}^{n}d_X(x_i,y_i)+\sum_{i=0}^{n-1}d_S(y_i,x_{i+1}).$$
	Define the \emph{semi-spherical metric} $d_{P_d}$ on $P_d(X)$ by
	$$d_{P_d}(x,y)=\inf\{l(\delta)\mid\text{ $\delta$ is a semi-simplicial path between $x$ and $y$}\}.$$
	One can show that $d_{P_d}$ is the largest metric such that $d_{P_d}(x,y)\leq d_S(x,y)$ and $d_{P_d}(z_1,z_2)\leq d_X(z_1,z_2)$ for any $x,y\in X(F)\subseteq P_d(X)$ for some finite $F\subseteq X$ and $z_1,z_2\in X$.
	
	The Rips complex becomes a metric space under the metric $d_{P_d}$. The canonical inclusion $i: X\to P_d(X)$ is a coarse equivalence for each $d\geq0$. One may refer to \cite[Chapter 7]{willett2020higher} for a proof. Then the coarse Baum-Connes conjecture states as follows:
	\begin{TEBCC}
		Let $X$ be a proper metric space. Then the coarse assembly map induced by the evaluation map is an isomorphism, i.e.,
		$$\mu_X:\lim_{d\to\infty}K_*(C^*_L(P_d(X)))\to\lim_{d\to\infty}K_*(C^*(P_d(X)))\cong K_*(C^*(X))$$
		is an isomorphism.
	\end{TEBCC}
	If we only consider the injectivity of the assembly map, we have the \emph{coarse Novikov conjecture}:
	\begin{TECNC}
		Let $X$ be a proper metric space. Then the coarse assembly map induced by the evaluation map is injective, i.e.,
		$$\mu_X:\lim_{d\to\infty}K_*(C^*_L(P_d(X)))\to\lim_{d\to\infty}K_*(C^*(P_d(X)))\cong K_*(C^*(X))$$
		is an injection.
	\end{TECNC}
	\section{A $C^*$-algebra associated with an admissible $\ell^p$-Busemann space}
	
	Fix $1\leq p<\infty$ and denote $B=\ell^p(\IN,\IR)$ to be the real Banach space and $\H=\ell^2(\IN,\IR)$ the real Hilbert space. The sign function $\operatorname{sgn}:\IR\to\IR$ is defined as
	$$\operatorname{sgn}(a)=\left\{\begin{aligned}&0&&,a=0,\\&a|a|^{-1}&&,a\ne0.\end{aligned}\right.$$
	Let $S(B)$ be the closed unit sphere of $B$. The Mazur map $\psi:S(B)\to S(\H)$, defined by
	$$\psi(a_1,\cdots,a_n,\cdots)=(\operatorname{sgn}(a_1)|a_1|^{p/2},\cdots,\operatorname{sgn}(a_n)|a_n|^{p/2},\cdots),$$
	is a uniform homeomorphism, i.e., the $r$-oscillation of $\psi$ and $\psi^{-1}$ tends to $0$ as $r$ tends to $0$, where the $r$-oscillation of $\psi$ is defined to be
	$$\omega_r(\psi)=\sup_{x,x'\in S(B),\|x-x'\|_p\leq r}\|\psi(x)-\psi(x')\|.$$
	The \emph{$p/2$-Hölder extension} of the Mazur map $\psi$ is a map $\phi:B\to\H$ defined by 
	$$\phi(x)=\left\{\begin{aligned}&0&&,x=0,\\&\|x\|^{p/2}_p\psi\left(\|x\|_p^{-1}x\right)&&,x\ne0.\end{aligned}\right.$$
	The extended Mazur map restricts to a uniform homeomorphism between bounded balls by \cite[Proposition 2.9]{OdellSchlumprecht1994} and Q. Cheng \cite[Theorem 2.4]{Cheng2016}:
	
	\begin{proposition}[\cite{OdellSchlumprecht1994,Cheng2016}]\label{extend Mazur uniform homeo}
		For any $R>0$, the extended Mazur map restricted on the ball $\phi|_{Ball_B(0,R)}:Ball_B(0,R)\to Ball_{\H}(0,R^{p/2})$ is a uniform homeomorphism, i.e.,
		$$\omega_r\left(\phi|_{Ball_B(0,R)}\right)\to 0\text{ as }r\to 0.$$
	\end{proposition}
	
	As a corollary, $\phi$ forms a homeomorphism between $B$ and $\H$. When $p=2$, the extended Mazur map is the identity map on $\H$.
	
	Let $\IR_+=[0,\infty)$, and let $\Cl(\H)$ be the complexified Clifford algebra of $\H$, one may refer to \cite[Section~3.1]{guo2024bott} for a definition of $\Cl(\H)$ for an infinite-dimensional Hilbert space. For any point $x$ in an $\ell^p$-Busemann space $M$, let $\ell^p_x$ be the smallest $\ell^p$ space into which the tangent cone $T_xM$ embeds. Note that the tip of the cone is at the origin. The Mazur map $\phi$ gives a homeomorphism from $\ell^p_x$ to a Hilbert space $\H_x$.
	We denote $C_0(\IR)$ by $\S$.
	\begin{definition}\label{defn:Pi-M}
		Given an admissible $\ell^p$-Busemann space $M$, we define the $*$-algebra
		\[
		\Pi(M) = \prod_{(x,t) \in M \times [0,\infty) } \Cl (\H_x \oplus t \IR ) \, ,
		\]
		where 
		\[
		t \IR = 
		\begin{cases}
			\IR \, , & t > 0 \\
			\{0\} \, , & t = 0 
		\end{cases}
		\, 
		\]
		and $\IR$ carries the canonical inner product independent of $t$. 
		We also define the $C^*$-algebra  
		\[
		\Pi_b(M) = \left\{ \sigma \in \Pi(M) \colon \sup_{(x,t) \in M \times [0,\infty)} \| \sigma({x,t}) \| < \infty \right\}
		\] 
		equipped with pointwise algebraic operations and the uniform norm. 
	\end{definition}
	
	\begin{definition}\label{defn_Cliffordmultiplier}
		Let $M$ be an admissible $\ell^p$-Busemann space. For any point $x_0 \in M$, we define the \emph{Clifford operator} $C_{x_0} \in \Pi (M)$ by
		\[
		C_{x_0}(x , t) = \left(\phi(-J_x(\log_{x}(x_0))), t \right) \in \H_{x} \times t \IR \subset \Cl({\H_{x}\oplus t \IR})
		\]
		for any $x \in M$ and $t \in [0,\infty)$, where $J_x: T_xM\to \ell^p_x$ is the bi-Lipschitz embedding, and $\phi$ is the Mazur map.
	\end{definition}
	
	Note that $C_{x_0}$ is an unbounded, self-adjoint operator. We denote by $C_0(\IR)_{ev}$ (or $C_0(\IR)_{odd}$, respectively) the subset of $C_0(\IR)$ of all even (or odd, respectively) functions. For any Hilbert space $H$, we next define the functional calculus for $\Cl(H)$: for any $f\in C_0(\IR)_{ev}$ and $v\in\H\subset\Cl(H)$, we define $f(v)=f(\|v\|)\in\IC\subset\Cl(H)$, and for $g\in C_0(\IR)_{odd}$, we define the functional calculus by
	$$g(v)=\left\{\begin{aligned}&0&&,v=0\\&g(\|v\|)\frac{v}{\|v\|}&&,v\ne0,\end{aligned}\right.$$
	then $g(v)$ is an element of $\Cl(H)$.
	\begin{definition}
		For any point $x_0\in M$, the \emph{Bott homomorphism}
		\[
		\beta_{x_0}: \S \to \Pi (M)
		\]
		is defined by the functional calculus
		\[
		\beta_{x_0}(f) (x,t) = f(C_{x_0}(x, t))
		\]
		for all $x \in M$, $t \in [0,\infty)$, and $f \in C_0(\IR)$. 
	\end{definition}
	Then $\beta_{x_0}$ is a graded $*$-homomorphism from $\S$ to the $*$-subalgebra $\Pi_b(M)$. Note that the Bott map takes $C_0(\IR)_\text{ev}$ into the subalgebra $C_0 (M \times [0,\infty))$ of the center of $\Pi_b(M)$ by definition of the functional calculus. We will write $\beta_{x_0}^M$ for $\beta_{x_0}$ if we want to emphasize the admissible $\ell^p$-Busemann space $M$. 
	
	\begin{definition}
		Let $M$ be an admissible $\ell^p$-Busemann space. The algebra $\mathcal{A}(M)$ is the $C^*$-subalgebra of $\Pi_b(M)$ generated by 
		\[
		\{ \beta_{x_0}(f) \colon x_0 \in M,\ f \in \S \ \} \, .
		\]
		We also define $ \mathcal{A}(M)_{\text{ev}} $ to be the $C^*$-subalgebra of $\mathcal{A}(M)$ generated by 
		\[
		\{ \beta_{x_0}(f) \colon x_0 \in M,\ f \in C_0(\IR)_\text{ev} \ \} \, . 
		\]
		Let $N$ be a closed convex subset of $M$. The algebra $\mathcal{A}(M,N)$ is defined to be the $C^*$-subalgebra of $\A(M)$ generated by
		$$\{\beta_x(f)\mid x\in N\subseteq M,f\in \S\}.$$
	\end{definition}
	
	\begin{lemma}\label{continuousBott}
		Fix $f\in\S$.\begin{enumerate}
			\item[(1)] For any $\varepsilon>0$, there exists $\delta>0$ such that
			$$\|\beta_x(f)-\beta_y(f)\|\leq\varepsilon$$
			whenever $x,y\in M$ satisfy $d(x,y)\leq\delta$.
			\item[(2)] As a corollary, let $\{x_n\}\in M$ be a sequence which converges to $x_0\in M$. Then
			$$\lim_{n\to\infty}\|\beta_{x_n}(f)-\beta_{x_0}(f)\|=0$$
			for each $f\in \S$.
		\end{enumerate}
	\end{lemma}
	
	\begin{proof}
		Let $J_x: T_xM\to \ell^p_x$ be the bi-Lipschitz embedding with Lipschitz constant $C>0$. Let $f_0$ and $f_1$ be the even and odd parts of $f$, respectively. It suffices to prove the claim for $f_0$ and $f_1$ separately. Since $f\in \S$, there exists $R>0$ such that $|f(t)|\leq\varepsilon/2$ for all $|t|>R$. 
		
		Even Case: By definition, for each $x,y,z\in M$, we have
		$$\|(\beta_x(f_0)-\beta_y(f_0))(z,t)\|=|f_0(\|C_x(z,t)\|)-f_0(\|C_y(z,t)\|)|.$$ If $(z,t)$ satisfies that both $\|C_x(z,t)\|$ and $\|C_y(z,t)\|$ are greater than $R$, then
		$$|f_0(\|C_x(z,t)\|)-f_0(\|C_y(z,t)\|)|\leq |f_0(\|C_x(z,t)\|)|+|f_0(\|C_y(z,t)\|)|\leq\varepsilon.$$
		Notice that $\|C_x(z,t)\|=(\|\phi(-J_z(\log_z(x)))\|^2+t^2)^{1/2}=(\|J_z(\log_z(x))\|^p_p+t^2)^{1/2}$.
		If one of $\|C_x(z,t)\|$ or $\|C_y(z,t)\|$ is less than $R$, without loss of generality, we assume that $\|C_x(z,t)\|<R$. Then $\|J_z(\log_z(x))\|^p_p< R^2$. Notice that
		\begin{equation*}\begin{split}
				\left|\|C_x(z,t)\|-\|C_y(z,t)\|\right|&=|\left(\|J_z(\log_z(x))\|^p_p+t^2\right)^{1/2}-\left(\|J_z(\log_z(y))\|^p_p+t^2\right)^{1/2}|\\
				&\leq \left(|\|J_z(\log_z(x))\|^p_p-\|J_z(\log_z(y))\|^p_p|\right)^{1/2},
		\end{split}\end{equation*}
		Recall that $\log_z$ is $1$-Lipschitz:
		$$ d_{T_zM}(\log_z(x),\log_z(y))\leq d_M(x,y).$$
		Since the tangent cone bi-Lipschitz embeds into $\ell^p_z$, we have
		$$\left|\|J_z(\log_z(x))\|_p-\|J_z(\log_z(y))\|_p\right|\leq \|J_z(\log_z(x))-J_z(\log_z(y))\|_p\leq Cd_{T_zM}(\log_z(x),\log_z(y))\leq Cd_M(x,y).$$
		Set $\delta_0>0$ satisfying that $|f_0(t_1)-f_0(t_2)|<\frac\varepsilon2$ for all $t_1,t_2$ with $|t_1-t_2|<\delta_0$. 
		As $y(t)=t^p$ is uniformly continuous on $[0,C^2R^{2/p}+1]$, there exists $\delta_1>0$ such that $|t_1^p-t_2^p|<\delta_0^2$ for all $t_1,t_2\in[0,C^2R^{2/p}+1]$ with $|t_1-t_2|<\delta_1$. Set $\delta_{ev}=\min\{\delta_1,1\}/C$. If $d_M(x,y)\leq \delta_{ev}$, then 
		$$\|J_z(\log_z(y))\|_p\leq Cd_M(z,y)\leq Cd_M(z,x)+Cd_M(x,y)\leq C^2\|J_z(\log_z(x))\|_p+Cd_M(x,y)\leq C^2R^{2/p}+1.$$
		Then both $\|J_z(\log_z(x))\|_p$ and $\|J_z(\log_z(y))\|_p$ are in $[0,C^2R^{2/p}+1]$, and we have
		$$|f_0(\|C_x(z,t)\|)-f_0(\|C_y(z,t)\|)|\leq\frac\varepsilon2.$$
		
		Odd Case: $tC_0(\IR)_{ev}$ is dense in $C_0(\IR)_{odd}$. We may assume that there exists $g\in C_0(\IR)_{ev}$ such that $f_1(t)=tg(t)$, then
		\begin{equation*}\begin{split}
				\|(\beta_x(f_1)-\beta_y(f_1))(z,t)\|=&\left\|g(\|C_x(z,t)\|)C_x(z,t)-g(\|C_y(z,t)\|)C_y(z,t)\right\|.
		\end{split}\end{equation*}
		If both $\|C_x(z,t)\|$ and $\|C_y(z,t)\|$ are greater than $R$, then
		$$\left\|g(\|C_x(z,t)\|)C_x(z,t)-g(\|C_y(z,t)\|)C_y(z,t)\right\|\leq |f_1(\|C_x(z,t)\|)|+|f_1(\|C_y(z,t)\|)|\leq\varepsilon.$$
		If one of $\|C_x(z,t)\|$ or $\|C_y(z,t)\|$ is less than $R$, we have	
		\begin{equation*}\begin{split}
				\|(\beta_x(f_1)-\beta_y(f_1))(z,t)\|\leq&|g(\|C_x(z,t)\|)-g(\|C_y(z,t)\|)|\left\|C_x(z,t)\right\|\\
				&+\|g(C_y(z,t))\|\cdot\|C_x(z,t)-C_y(z,t)\|.
		\end{split}\end{equation*} 
		With a similar argument as above, we can find a $\delta_3>0$ such that $|g(\|C_x(z,t)\|)-g(\|C_y(z,t)\|)|\leq \frac \varepsilon{4R}$. For the second term, we have
		$$\|C_x(z,t)-C_y(z,t)\|=\|\phi(-J_z(\log_z(x)))-\phi(-J_z(\log_z(y)))\|.$$
		By Proposition \ref{extend Mazur uniform homeo}, there exists $0<\delta_4<1$ associated with $\phi|_{Ball_B(0,C^2R^{2/p}+1)}$ such that if $d_M(x,y)\leq \delta_4$, then $\|\phi(-J_z(\log_z(x)))-\phi(-J_z(\log_z(y)))\|\leq\frac\varepsilon{4\|g\|}$. Set $\delta_{odd}=\min\{\delta_3,\delta_4\}$, we conclude that if $d_M(x,y)\leq\delta_{odd}$, then
		$$\|(\beta_x(f_1)-\beta_y(f_1))(z,t)\|\leq\frac\varepsilon2.$$
		Set $\delta=\min\{\delta_{ev},\delta_{odd}\}$, then
		$$\|\beta_x(f)-\beta_y(f)\|\leq\|\beta_x(f_0)-\beta_y(f_0)\|+\|\beta_x(f_1)-\beta_y(f_1)\|\leq\varepsilon.$$
	\end{proof}
	Since $M$ is a geodesic space, given any $x_0,x_1\in M$, there is a geodesic segment $(x_s)_{s\in [0,1]}$ in $M$ connecting $x_0$ and $x_1$. This gives the following corollary: 
	\begin{corollary}\label{bott-homotopy}
		For any two points $x_0, x_1 \in M$, the Bott homomorphisms
		\[
		\beta_{x_0}, \beta_{x_1} \colon \S \to \A(M)
		\]
		are homotopic to each other. 
	\end{corollary}
	\begin{corollary}\label{directlimit}
		Let $\{M_n\}_{n\in\IN}$ be an increasing sequence of closed convex subsets of $M$ with $\bigcup_{n\in\IN}M_n$ dense in $M$. Then
		\begin{enumerate}
			\item[(1)] $\A(M,M_n)\subseteq\A(M,M_{n+1})$,
			\item[(2)] $\lim_{n\to\infty}\A(M,M_n)=\A(M,\bigcup_{n\in\IN}M_n)=\A(M,\overline{\bigcup_{n\in\IN}M_n})$.
		\end{enumerate}
	\end{corollary}
	
	\begin{proposition}
		Let $\widehat{\mathcal{A}(M)_{\text{ev}}}$ be the spectrum of the central subalgebra $\mathcal{A}(M)_{\text{ev}}$. Then $\mathcal{A}(M)$ is an $\widehat{\mathcal{A}(M)_{\text{ev}}}$-{$C^*$-algebra}. 
	\end{proposition}
	\begin{proof}
		We already know that $\mathcal{A}(M)_{\text{ev}}$ is a sub-$C^*$-algebra of the center $Z(\mathcal{A}(M))$. It suffices to show that $\mathcal{A}(M)_{\text{ev}} \cdot \mathcal{A}(M)$ is dense in $\mathcal{A}(M)$, which follows from the definition of $\mathcal{A}(M)$ and the fact that every $f\in\S$ can be approximated as a product $f = f_1 f_2$ where $f_1 \in \S_\text{ev}$ and $f_2 \in \S$.
	\end{proof}
	It is clear that $\widehat{\mathcal{A}(M)_{\text{ev}}}=M\times [0,\infty)$, we equip $M\times [0,\infty)$ with the weak-$*$ topology so that it is a locally compact Hausdorff space. If $M$ is finite-dimensional, then the weak-$*$ topology on $M\times [0,\infty)$ coincides with the product topology.
	
	\begin{lemma}\label{Aev=c0}
		Let $M$ be a complete, connected, and simply connected finite-dimensional $\ell^p$-Finsler manifold. Then $\A_{\text{ev}} (M)$ coincides with $C_0(M \times [0,\infty))$ as $C^*$-subalgebras of $\Pi_{b} (M)$.
	\end{lemma}
	
	\begin{proof}
		It is clear that $\A_{\text{ev}} (M) \subset C_0(M \times [0,\infty))$. Let $f\in\S_{\text{ev}}$, then
		$$\beta_x(f)(z,t)=f_0(\|C_x(z,t)\|)=f_0((\|\phi(-J_z(\log_z(x)))\|^2+t^2)^{1/2})=f_0((\|J_z(\log_z(x))\|_p^p+t^2)^{1/2}).$$ 
		For any two distinct points $(x_1, t_1)$ and $(x_2 , t_2)$ in $M \times [0,\infty)$ with $t_1 \leq t_2$, we have
		\[
		d(x_1,x_1)^2 + t_1^2 = t_1^2 < d(x_1,x_2)^2 + t_2^2.
		\]
		We can choose $f \in \S_\text{ev}$ so that 
		\[
		\beta_{x_1} (f) (x_1, t_1) \not= \beta_{x_1} (f) (x_2, t_2) \, .
		\]
		So the generating set 
		\[
		\left\{ \beta_{x_0} (f) \colon f \in \S_\text{ev}, x_0 \in M \right\}
		\]
		separates points. Hence the Stone-Weierstrass theorem implies $\A_{\text{ev}}(M) = C_0(M \times [0,\infty))$. 
	\end{proof}
	
	\begin{proposition}
		Let $N$ be a finite-dimensional $\ell^p$ space, viewed as an $\ell^p$-Busemann space with each tangent cone identified canonically with $N$ itself, and the logarithm map being the identity. Assume $N$ is homeomorphic to the Hilbert space $H$ via the Mazur map. Then the $C^*$-subalgebra $\A(N)$ coincides with 
		\begin{equation}
			\left\{ f \in C_0 \big( N \times [0, \infty) , \Cl (H \oplus \IR ) \big) \colon f( x , 0 ) \in \Cl (H) \text{ for all } x \in N \right\} \, ,
		\end{equation}
		where the embedding $\Cl(H) \hookrightarrow \Cl (H\oplus \IR )$ is induced by embedding $H$ into the first factor of $H \oplus \IR$.  
	\end{proposition}
	
	\begin{proof}
		Denote the $C^*$-algebra in the equation by $\mathcal{B}$. By Lemma~\ref{Aev=c0}, $\A(N)$ and $\mathcal{B}$ are both $( N \times [0, \infty) )$-$C^*$-algebras. 
		Now for $(x, t) \in N \times [0,\infty)$, we see that  
		\[
		\left\{ C_{x_0} (x,t) \colon x_0 \in N \right\} = H \times \{ t \} \subset H \oplus t \IR \subset \Cl (H \oplus t \IR) \, .
		\]
		As $( N \times [0, \infty) )$-$C^*$-algebras, the fibers of $\A(N)$ and $\mathcal{B}$ at $(x,t)$ both coincide with $\Cl (H\oplus t \IR )$. Therefore $\A(N) = \mathcal{B}$ as $C^*$-subalgebras of $\Pi(N)$. 
	\end{proof}
	
	It follows that for an $\ell^p$ space $N$, viewed as an $\ell^p$-Busemann space, our $C^*$-algebra $\A(N)$ is the same as the one defined in \cite{guo2024bott}. For a Hilbert-Hadamard space $M$ with the Mazur map being the identity, our $\A(M)$ also coincides with the one defined in \cite{gong2021novikov}. By \cite[Corollary 3.13]{guo2024bott}, we conclude that:
	
	\begin{corollary}Let $N$ be a finite-dimensional $\ell^p$ space, then the Bott map induces an isomorphism on $K$-theory:
		$$(\beta_x^N)_*:K_*(\S)\xrightarrow{\cong} K_*(\A(N))$$
		for any $x\in N$.
	\end{corollary}
	Let $M$ be a complete, simply connected $\ell^p$-Finsler-Busemann manifold, and let $V$ be the tangent space at $x_0\in M$, viewed as a finite-dimensional normed vector space with the Finsler norm $F_{\epsilon}$. The Finsler-Cartan-Hadamard theorem gives a diffeomorphism from $V$ to $M$. We aim to establish an isomorphism between $\A(M)$ and $\A(V)$ in analogy with \cite[Proposition 5.21]{gong2021novikov}. However, there is one important difference: in the Hilbert-Hadamard space case, the authors use the polar decomposition of $D_v\exp_{x_0}$ to obtain a continuous family of linear isometries between tangent spaces. But for our Berwald $\ell^p$-Finsler manifold, the natural analogue is to use parallel transport along radial geodesics, which preserves the Finsler norm.
	
	\begin{proposition}\label{prop_AofM-fin-dim}
		Let $M$ be a complete, simply connected $\ell^p$-Finsler-Busemann manifold, and let $V$ be the tangent space at $x_0\in M$. Then there is an isomorphism 
		\[
		\Phi \colon \A(V) \to \A(M)
		\]
		which intertwines the Bott homomorphisms in the sense that $\Phi \circ \beta_0^{V} = \beta_{x_0}^{M}$. 
	\end{proposition}

	\begin{proof}
		Let $V=T_{x_0}M$. By the Finsler-Cartan-Hadamard theorem, $\exp_{x_0}:V\to M$ is a diffeomorphism. For each $v\in V$, let $P_v:V\to T_{\exp_{x_0}(v)}M$ be parallel transport along the radial geodesic from $x_0$ to $\exp_{x_0}(v)$. Since $M$ is Berwald, parallel transport preserves the Finsler norm, so $P_v$ is a linear isometry. Thus
		\[
		F_{\varepsilon,\exp_{x_0}(v)}(P_vw)
		=F_{\varepsilon,x_0}(w)
		\]
		for every $w\in V$, and
		$$
		P_v(v)=-\log_{\exp_{x_0}(v)}(x_0).
		$$
		Indeed, the velocity of $\gamma_v$ at its endpoint is $P_v(v)$, while the geodesic from $\exp_{x_0}(v)$ back to $x_0$ has initial velocity $-P_v(v)$. We first construct a compatible family of embeddings of the tangent spaces into $\ell_n^p$. Let
		\[
		J_0: (V,F_{\epsilon})\longrightarrow\ell_n^p
		\]
		be the identity map. There exist constants $C>0$ such that
		\[
		C^{-1}F_{\varepsilon,x_0}(w)
		\leq \|J_0w\|_p
		\leq CF_{\varepsilon,x_0}(w)
		\]
		for all $w\in V$. For $x=\exp_{x_0}(v)$, define
		\[
		J_x=J_0\circ P_v^{-1}:
		T_xM\longrightarrow\ell_n^p.
		\]
		Then
		\[
		\begin{split}
			C^{-1}F_{\varepsilon,x}(w)
			&=C^{-1}F_{\varepsilon,x_0}(P_v^{-1}w)\\
			&\leq \|J_0P_v^{-1}w\|_p\\
			&\leq CF_{\varepsilon,x_0}(P_v^{-1}w)\\
			&=CF_{\varepsilon,x}(w).
		\end{split}
		\]
		Hence the maps $J_x$ are homogeneous and bi-Lipschitz, with Lipschitz constant independent of $x$. Moreover, we have
		\[
		J_{\exp_{x_0}(v)}P_v=J_0.
		\]
		The family $J_x$ depends continuously on $x$ since parallel transport along radial geodesics depends continuously on $v$. Let $\phi:\ell_n^p\to H=\ell_n^2$ be the Mazur map. We use $J_x$ and $\phi$ to identify the Hilbert space $H_x$ associated to $T_xM$ with $H$.
		Define $U_v:H_V\to H_{\exp_{x_0}(v)}$ by
		$$
		U_v=\phi_{\exp_{x_0}(v)}J_{\exp_{x_0}(v)} P_v J_0^{-1} \phi_V^{-1}.
		$$
		By construction, the induced map between the Hilbert spaces associated with $x_0$ and $\exp_{x_0}(v)$ is
		\[
		\begin{split}
			U_v
			&=\phi\circ J_{\exp_{x_0}(v)}
			\circ P_v\circ J_0^{-1}\circ\phi^{-1}\\
			&=\phi\circ J_0\circ P_v^{-1}P_v
			\circ J_0^{-1}\circ\phi^{-1}\\
			&=\operatorname{id}_H,
		\end{split}
		\]
		We now construct the desired isomorphism. For $x=\exp_{x_0}(v)$, use the unitary $U_v$ to identify the Clifford fibers. Define
		\[
		\begin{split}
			\widetilde{\Phi}:\Pi(V)&\longrightarrow \Pi(M),\\
			\widetilde{\Phi}(\sigma)(x,t)
			&=
			\operatorname{Cl}(U_{\log_{x_0}(x)}\oplus\operatorname{id}_{\mathbb R})
			\big(\sigma(\log_{x_0}(x),t)\big),
			\qquad v=\log_{x_0}(x).
		\end{split}
		\]
		Here, $\operatorname{Cl}(U_{\log_{x_0}(x)}\oplus\operatorname{id}_{\mathbb R})$ denotes the isomorphism between Clifford algebras induced by $U_{\log_{x_0}(x)}\oplus\operatorname{id}_{\mathbb R}$. It remains to verify that $\widetilde{\Phi}$ intertwines the Bott homomorphisms based at $0$ and $x_0$. We have
		\[
		\begin{split}
			J_x(-\log_x(x_0))
			&=J_xP_v(v)\\
			&=J_0(v).
		\end{split}
		\]
		It follows that
		\[
		\begin{split}
			C_{x_0}^M(\exp_{x_0}(v),t)
			&=\phi(J_x(-\log_x(x_0)))\oplus te\\
			&=\phi(J_0(v))\oplus te\\
			&=C_0^V(v,t),
		\end{split}
		\]
		i.e., $C_{x_0}^M(x,t)=C_0^V(\log_{x_0}(x),t)$.
		Applying the continuous functional calculus, we get
		\[
		\widetilde{\Phi}(\beta_0^V(f))
		=\beta_{x_0}^M(f)
		\]
		for every $f\in\mathcal S$.
		Finally, $\widetilde\Phi$ identifies the central even subalgebras through $\exp_{x_0}$, and the fibers of both $\widetilde\Phi(\mathcal A(V))$ and $\mathcal A(M)$ at $(x,t)$ are $\Cl(H_V\oplus t\mathbb R)$.
		Consequently,
		$$
		\widetilde\Phi(\mathcal A(V))=\mathcal A(M).
		$$
		The restriction $\Phi=\widetilde\Phi|_{\mathcal A(V)}$ is the required $*$-isomorphism and satisfies $\Phi\circ\beta_0^V=\beta_{x_0}^M$. Hence
		\[
		\Phi\circ\beta_0^V=\beta_{x_0}^M,
		\]
		as required.
	\end{proof}
	\begin{proposition}
		Let $V=\mathbb{R}^n$ be a finite-dimensional real vector space, and let $1< p<\infty$. Equip $V$ with either the $\ell^p$-norm
		\[
		\|y\|_p=\left(\sum_{i=1}^n |y_i|^p\right)^{1/p}
		\]
		or the norm
		\[
		F_\varepsilon(y)=
		\left(
		\sum_{i=1}^n
		(y_i^2+\varepsilon^2|y|_2^2)^{p/2}
		\right)^{1/p}.
		\]
		Then, using the homogeneous bi-Lipschitz identity embedding $J: (V,F_{\varepsilon}) \to (V,\ell^p)$, the associated $C^*$-algebras coincide:
		\[
		\mathcal A(V,\ell^p)=\mathcal A(V,F_\varepsilon).
		\]
	\end{proposition}
	
	\begin{proof}
		We identify $V$ with $\mathbb R^n$, and let $\phi:\ell_n^p\to H$ be the extended Mazur map appearing in the definition of the associated $C^*$-algebras, with
		\[
		J:V\longrightarrow \ell_n^p,\qquad
		J(y_1,\ldots,y_n)=(y_1,\ldots,y_n).
		\]
		The geodesics for either norm are straight lines.
		In particular, for every $x,y\in V$, we have
		\[
		\log_x^{\ell^p}(y)=y-x
		\quad\text{and}\quad
		\log_x^{F_\varepsilon}(y)=y-x.
		\]
		Thus
		\[
		-\log_x^{\ell^p}(y)
		=
		-\log_x^{F_\varepsilon}(y)
		=x-y.
		\]
		Recall that the Clifford operator associated with $y\in V$ is defined by
		\[
		C_y(x,t)
		=
		\bigl(\phi(-J(\log_x(y))),t\bigr)
		\in H\oplus\mathbb R.
		\]
		Therefore, for every $(x,t)\in V\times[0,\infty)$,
		\[
		\begin{aligned}
			C_y^{F_\varepsilon}(x,t)
			&=
			\bigl(\phi(-J(\log_x^{F_\varepsilon}(y))),t\bigr)\\
			&=
			\bigl(\phi(J(x-y)),t\bigr)\\
			&=
			\bigl(\phi(-J(-\log_x^{\ell^p}(y))),t\bigr)\\
			&=
			C_y^{\ell^p}(x,t).
		\end{aligned}
		\]
		It follows immediately from the functional calculus that, for every $a\in V$ and $f\in C_0(\mathbb R)$,
		\[
		\beta_a^{F_\varepsilon}(f)
		=
		\beta_a^{\ell^p}(f)
		\]
		Finally, we have
		\[
		\begin{aligned}
			\mathcal A(V,F_\varepsilon)
			&=\left\{
			\beta_a^{F_\varepsilon}(f):
			a\in V,\ f\in C_0(\mathbb R)
			\right\}\\
			&=\left\{
			\beta_a^{\ell^p}(f):
			a\in V,\ f\in C_0(\mathbb R)
			\right\}\\
			&=\mathcal A(V,\ell^p).
		\end{aligned}
		\]
		This completes the proof.
	\end{proof}
	\begin{corollary}\label{bottisoM}Let $M$ be a finite-dimensional $\ell^p$-Finsler-Busemann manifold, then the Bott map induces an isomorphism on $K$-theory:
		$$(\beta_x^M)_*:K_*(\S)\xrightarrow{\cong} K_*(\A(M))$$
		for any $x\in M$.
	\end{corollary}
	
	Let $M$ be an $\ell^p$-Busemann space, and $N$ be a closed convex subset of $M$. Next, we compare the algebras $\A(M)$ and $\A(M,N)$. Let
	\[
	\Pi(M,N) = \left\{ \sigma \in \prod_{(x,t) \in N \times [0,\infty) } \Cl(H_x M \oplus t \IR ) \colon \sup_{(x,t) \in N \times [0,\infty)} \| \sigma({x,t}) \| < \infty \right\}.
	\]
	There is a natural quotient map 
	\[
	\pi_{M,N} \colon \Pi_b(M) \to \Pi_b(M,N)
	\]
	and a natural embedding
	\[
	\iota_{M, N} \colon \Pi_b(N) \to \Pi_b(M,N)
	\]
	induced by the embeddings $\H_x N \hookrightarrow \H_x M$ for all $x \in N$. 
	
	\begin{lemma}\label{AMN-AM}
		Let $N$ be a closed convex subset of $M$. Then for any $x_0 \in N$, we have
		\[
		\pi_{M,N} \circ \beta_{x_0}^M= \iota_{M,N} \circ \beta_{x_0}^N.
		\]
		In particular,
		\[
		\pi_{M,N}(\A(M,N)) = \iota_{M,N}(\A(N)).
		\]
	\end{lemma}
	
	\begin{proof}
		We have
		\begin{align*}
			\pi_{M,N} \circ \beta_{x_0}^Mf(x,t) = f\left(C^M_{x_0}(x,t) \right)
		\end{align*}
		and 
		\begin{align*}
			\iota_{M,N} \circ \beta_{x_0}^Nf(x,t)=f\left(C^N_{x_0}(x,t) \right).
		\end{align*}
		Since $N$ is convex and closed, the geodesic segment in $M$ connecting $x_0$ to $x$ coincides with the geodesic segment in $N$ connecting $x_0$ to $x$. So they give the same element.
	\end{proof}
	
	\begin{proposition}\label{bottinjective}
		Suppose that $M$ is an admissible $\ell^p$-Busemann space. Let $M_1 \subset M_2 \subset \ldots$ be closed convex subsets isometric to finite-dimensional $\ell^p$-Finsler-Busemann manifolds, whose union is dense in $M$. Then for any $x_0 \in M$, the Bott homomorphism 
		\[
		\beta_{x_0}  \colon \S \to \A(M)
		\]
		induces injections on $K$-theory, i.e., 
		\[
		K_i (\S) \xlongrightarrow{(\beta_{x_0} )_*} K_i (\A(M)) 
		\] 
		is injective. 
	\end{proposition}
	\begin{proof}
		By Corollary~\ref{bott-homotopy}, we may assume without loss of generality that $x_0 \in M_1$. By Corollary~\ref{directlimit} and the continuity of the $K$-theory functor with regard to direct limits, it suffices to show that  
		\[
		\beta_{x_0} \colon \S \to \A(M,M_k)
		\]
		induces an injection on $K$-theory for every $k>0$. By Lemma~\ref{AMN-AM} we have a commutative diagram:
		\[
		\xymatrix{
			\S \ar[r]^{\beta_{x_0} } \ar[d]_{\beta_{x_0}^{M_k} } & \A(M,M_k)\ar[d]^{\pi_{M,M_k}} \\
			\A(M_k)\ar[r]_-{\iota_{M,M_k}} &	 \iota_{M,M_k}\left( \A(M_k)\right) = \pi_{M,M_k} \left(\A(M,M_k)\right)
		}
		\]
		The lower horizontal arrow is an isomorphism by Lemma~\ref{AMN-AM}, and the left vertical arrow induces isomorphisms on $K$-theory by Corollary~\ref{bottisoM}. It follows that $\beta_{x_0}$ induces injections on $K$-theory. 
	\end{proof}
	\section{The twisted algebras and the Bott homomorphisms}
	\subsection{The twisted algebras}
	Throughout the rest of this paper, we assume that $X$ is a proper metric space with bounded geometry which coarsely embeds into an admissible $\ell^p$-Busemann space $M$. Let $f:X\to M$ be the coarse embedding and $N\subset X$ be a countable dense subset of $X$. 
	\begin{definition}
		Let $X$ be a proper metric space. A \emph{(Roe algebraic) coefficient system} for $X$ is a triple $(N,M,f)$ consisting of the following data:
		\begin{itemize}
			\item[(1)] a countable dense subset $N\subseteq X$,
			\item[(2)] an admissible $\ell^p$-Busemann space $M$,
			\item[(3)] a coarse embedding $f:N\to M$.
		\end{itemize}
	\end{definition}
	
	\begin{definition}\label{twistedroe}
		The \emph{algebraic twisted Roe algebra} for the coefficient system $(N,M,f)$, denoted by $\mathbb C[X,\mathcal A(M)]$, is the set of all bounded functions $T$ on $N\times N$ such that:\begin{itemize}
			\item[(1)] $T(x,y)\in\K(\H_0)\otimes\mathcal A(M)$ for each $x,y\in N$ and there exists $R>0$ such that
			$$\textup{supp}(T(x,y))\subseteq B(f(x),R)\subseteq M\times\mathbb R_+,$$
			where $B(v,R)=\{(w,t)\in M\times \mathbb{R}_{+}\mid d_{M}(w,v)^{2}+t^{2}<R^{2}\}$,
			\item[(2)] there exists $r>0$ such that $T(x,y)=0$ for all $x,y\in N$ with $d(x,y)>r$,
			\item[(3)] for each bounded set $B\subseteq X$, the set
			$$\{(x,y)\in (N\times N)\cap (B\times B)\mid T(x,y)\ne0\}$$
			is finite,
			\item[(4)] there exists $L>0$ such that
			$$\#\{y\in N\mid T(x,y)\ne0\}\leq L,\quad\#\{y\in N\mid T(y,x)\ne0\}\leq L$$
			for all $x\in N$.
	\end{itemize}\end{definition}
	The algebraic structure on $\mathbb C[X,\mathcal A(M)]$ is given by
	$$(TS)(x,y)=\sum_{z\in N}T(x,z)S(z,y)$$
	$$(T^*)(x,y)=(T(y,x))^*.$$
	The algebraic twisted Roe algebra is represented as adjointable operators on the following Hilbert $C^*$-module over $\mathcal A(M)$: $$E_{\mathcal A(M)}=\ell^2(N)\otimes\H_0\otimes\mathcal A(M)$$
	$$T(\delta_x\otimes v\otimes a)=\sum_{y\in N}\delta_y\otimes T(y,x)(v\otimes a),$$
	
	\begin{definition}
		The \emph{twisted Roe algebra} for the coefficient system $(N,M,f)$, denoted by $C^*(X,\mathcal A(M))$, is the norm closure of $\mathbb C[X,\mathcal A(M)]$ in the $C^*$-algebra of all adjointable module homomorphisms.
	\end{definition}
	
	\begin{definition}\label{twistedlocalization}
		Define $\mathbb C_L[X,\mathcal A(M)]$ to be the set of all bounded, uniformly norm-continuous functions
		$$g:\mathbb R_+\to\mathbb C[X,\mathcal A(M)]$$
		such that\begin{itemize}
			\item[(1)] there exists a bounded function $r:\mathbb R_+\to\mathbb R_+$ with $\lim_{s\to\infty}r(s)=0$ such that $(g(s))(x,y)=0$ whenever $d(x,y)>r(s)$,
			\item[(2)] there exists $R>0$ such that $\textup{supp}((g(s))(x,y))\subseteq B(f(x),R)\subseteq M\times\mathbb R_+$ for all $s\in\mathbb R_+$.
	\end{itemize}\end{definition}
	
	\begin{definition}
		The \emph{twisted localization algebra} for the coefficient system $(N,M,f)$, denoted by $C^*_L(X,\mathcal A(M))$, is the completion of $\mathbb C_L[X,\mathcal A(M)]$ with respect to the norm
		$$\|g\|=\sup_{s\in\mathbb R_+}\|g(s)\|.$$
	\end{definition}
	
	Let $P_d(X)$ be the Rips complex of $X$. We can choose a countable dense subset $N_d\subset P_d(X)$ such that $N\subset N_d\subset N_{d'}$ for $d<d'$. Since there is a coarse equivalence $N_d\xrightarrow{i_d}X$, $(N_d,M,f\circ i_d)$ is a coefficient system for $P_d(X)$. We can then define the twisted Roe algebra for $P_d(X)$, denoted by $C^*(P_d(X),\mathcal A(M))$, and the twisted localization algebra for $P_d(X)$, denoted by $C^*_L(P_d(X),\mathcal A(M))$, respectively, as before.
	\subsection{The Bott maps}
	For any $g\in\S$ and $t\in[1,\infty)$, we define $g_t$ to be the rescaling function given by
	$$(g_t)(s)=g(t^{-1}s).$$ 
	\begin{definition}
		For each $d\geq 0$ and $t\in[1,\infty)$, define a map
		$$\beta_t:\S\otimes\IC[P_d(X)]\to\IC[P_d(X),\A(M)]$$
		for each $g\in\S$ and $T\in\IC[P_d(X)]$ by the formula
		$$(\beta_t(g\otimes T))(x,y)=T(x,y)\otimes(\beta(x))(g_t),$$
		where
		$$(\beta(x))(g_t)=\beta_{f(x)}^M(g_t).$$
	\end{definition}
	
	\begin{definition}
		For each $d\geq 0$ and $t\in[1,\infty)$, define a map
		$$(\beta_L)_t:\S\otimes\IC_{L}[P_d(X)]\to\IC_{L}[P_d(X),\A(M)]$$
		for each $g\in\S$ and $f\in\IC_{L}[P_d(X)]$ by the formula
		$$((\beta_L)_t(g\otimes f))(r)=\beta_t(g\otimes f(r)),$$
		for each $r\in\IR_+$.
	\end{definition}
	
	The following lemma is proved similarly to \cite[Lemma 4.9]{guo2024hilbert}.
	\begin{lemma}\label{bottasym}
		The maps $\beta_t$ and $(\beta_L)_t$ extend to asymptotic morphisms as $t\rightarrow \infty$:
		$$\beta:\S\otimes C^*(P_d(X))\to C^*(P_d(X),\A(M)),$$
		$$\beta_L:\S\otimes C^*_L(P_d(X))\to C^*_L(P_d(X),\A(M)).$$
	\end{lemma}
	\begin{proof}
		Below we show that
		$$\|\beta_t(TS)-\beta_t(T)\beta_t(S)\|\to 0$$
		as $t\to\infty$ for any $T,S\in\S\otimes\mathbb C_{\rm f}[P_d(X)]$.
		\begin{equation*}
			\begin{aligned}
				&\beta_{t}((g_{1}\otimes T_{1})(g_{2}\otimes T_{2}))(x,x')-(\beta_{t}(g_{1}\otimes T_{1})\cdot\beta_{t}(g_{2}\otimes T_{2}))(x,x')\\
				&=(T_{1}T_{2})(x,x')\otimes \beta(x)(g_{1}g_{2})_{t} -\sum_{\xi\in N}
				(\beta_{t}(g_{1}\otimes T_{1})(x,\xi))\cdot (\beta_{t}(g_{2}\otimes T_{2}))(\xi,x')\\
				&=\sum_{\xi\in N}T_{1}(x,\xi)T_{2}(\xi,x')\otimes (\beta(x)(g_{1})_{t})(\beta(x)(g_{2})_{t})\\
				&\quad-\sum_{\xi\in N}T_{1}(x,\xi)\otimes (\beta(x)(g_{1})_{t}) \cdot T_{2}(\xi,x')
				\otimes (\beta(\xi))(g_{2})_{t}\\
				&=\sum_{\xi\in N}T_{1}(x,\xi)T_{2}(\xi,x')\otimes \beta(x)(g_{1})_{t}
				[\beta(x)-\beta(\xi)](g_{2})_{t}
			\end{aligned}
		\end{equation*}
		Since $X$ has bounded geometry, by \cite[Lemma~3.4]{gong2008geometrization}, it suffices to show that for any $r>0$ and $\varepsilon>0$, there exists $K>0$ such that
		$$\|\beta(x)(g_t)-\beta(y)(g_t)\|\leq \varepsilon$$
		for all $t>K$ and $x,y\in N_d$ with $d(x,y)<r$. Since $f: X\to M$ is a coarse embedding, there exists $r'>0$ such that $d(f(x),f(y))<r'$ for all $x,y$ with $d(x,y)<r$.
		For any $s\in[0,1]$, one has
		\[\begin{split}&\left\|\left(\beta(x)(g_t)-\beta(y)(g_t)\right)\right\|=\|\beta_{f(x)}^M(g_t)-\beta^{M}_{f(y)}(g_t)\|\\
			=&\sup_{(z,r)\in M\times\mathbb R_+}\left\|g_t\left(C^{M}_{f(x)}(z,r)\right)-g_t\left(C^{M}_{f(y)}(z,r)\right)\right\|\\
			=&\sup_{(z,r)\in M\times\mathbb R_+}\left\|g\left(\frac{\phi(-J_z(\log_z(f(x))))}t,\frac rt\right)-g\left(\frac{{\phi(-J_z(\log_z(f(y))))}}t,\frac rt\right)\right\|\\
			=&\sup_{(z,r)\in M\times\mathbb R_+}\left\|g\left(\phi(-t^{-2/p}J_z(\log_z(f(x)))),\frac rt\right)-g\left({\phi(-t^{-2/p}J_z(\log_z(f(y))))},\frac rt\right)\right\|.\\
		\end{split}\]
		Let $\gamma_1$ be the geodesic segment from $z$ to $f(x)$, and $\gamma_2$ be the geodesic segment from $z$ to $f(y)$. We have
		$$t^{-2/p}\log_z(f(x))=[\gamma_1, t^{-2/p}d(z,f(x))]\in T_zM$$ 
		and 
		$$t^{-2/p}\log_z(f(y))=[\gamma_2, t^{-2/p}d(z,f(y))]\in T_zM.$$
		
		Let $x'=\gamma_1(t^{-2/p}d(z,f(x)))$ and $y'=\gamma_2(t^{-2/p}d(z,f(y)))$. Note that $J_z$ is homogeneous, so we have 
		\[\begin{split}\left\|\left(\beta(x)(g_t)-\beta(y)(g_t)\right)\right\|&=\sup_{(z,r)\in M\times\mathbb R_+}\left\|g\left(\phi(-J_z(\log_z(x'))),\frac rt\right)-g\left({\phi(-J_z(\log_z(y'))),\frac rt}\right)\right\|\\
			&=\left\|\left(\beta(x')(g)-\beta(y')(g)\right)\right\|.
		\end{split}\]
		Note that by Busemann convexity, $d(x',y')\leq t^{-2/p}d(f(x),f(y))$. Using the fact $\lim\limits_{t\to\infty}t^{-2/p}=0$, that $J_z$ is bi-Lipschitz, and Lemma \ref{continuousBott}, we conclude that $\left\|\left(\beta(x)(g_t)-\beta(y)(g_t)\right)\right\|\to 0.$
	\end{proof}
	
	The asymptotic morphisms
	$$\beta_t:\S\otimes C^*(P_d(X))\to C^*(P_d(X),\A(M))$$
	$$(\beta_L)_t:\S\otimes C^*_{L}(P_d(X))\to C^*_{L}(P_d(X),\A(M))$$
	induce homomorphisms on $K$-theory:
	$$\beta_*:K_*(\S\otimes C^*(P_d(X)))\to K_*(C^*(P_d(X),\A(M))),$$
	$$(\beta_L)_*:K_*(\S\otimes C^*_{L}(P_d(X)))\to K_*(C^*_{L}(P_d(X),\A(M))).$$
	\begin{theorem}\label{BottL}
		Let $M$ be a finite-dimensional $\ell^p$-Finsler-Busemann manifold, then for any $d> 0$, the Bott map
		$$(\beta_L)_*:K_*(\S\otimes C^*_{L}(P_d(X)))\to K_*(C^*_{L}(P_d(X),\A(M)))$$
		is an isomorphism.
	\end{theorem}
	\begin{proof}
		The $K$-theory of the localization algebra is invariant under strong Lipschitz homotopy equivalence (see \cite{yu1997localization} and \cite{yu2000coarse}). By a Mayer-Vietoris sequence argument and induction on the dimension of the skeletons, the general case can be reduced to the zero-dimensional case. Let $X^{(0)}\subseteq P_d(X)$ be the $0$-skeleton of $P_d(X)$, we prove
		$$(\beta_L)_*:K_*\left(\S\otimes C^*_{L}\left(X^{(0)}\right)\right)\to K_*\left(C^*_{L}\left(X^{(0)},\A(M)\right)\right)$$
		is an isomorphism.
		
		Notice that
		$$K_*\left(\S\otimes C^*_{L}\left(X^{(0)}\right)\right)=\prod_{x\in X^{(0)}}K_*(\S\otimes C^*_L(\{x\}))$$
		$$K_*\left(C^*_{L}\left( X^{(0)},\A(M)\right)\right)=\prod_{x\in X^{(0)}}K_*(C^*_L(\{x\},\A(M))).$$
		Moreover, one can see that $K_*(\S\otimes C^*_L(\{x\}))\cong K_*(\S\otimes\K)$ and $K_*(C^*_L(\{x\},\A(M)))\cong K_*(\K\otimes\A(M))$, and the Bott map $(\beta_L)_*$ coincides with $\beta_{f(x)}$ when $(\beta_L)_*$ restricts on $K_*(\{x\})\cong K_*(\S)$ to $K_*(C^*_L(\{x\},\A(M)))\cong K_*(\A(M))$. Then the theorem follows from Corollary \ref{bottisoM}.
	\end{proof}
	
	\begin{theorem}\label{BLinjective}
		Let $M$ be an admissible $\ell^p$-Busemann space, then for any $d> 0$, the Bott map
		$$(\beta_L)_*:K_*(\S\otimes C^*_{L}(P_d(X)))\to K_*(C^*_{L}(P_d(X),\A(M)))$$
		is a $\mathbb{Q}$-injection.
	\end{theorem}
	The proof depends on the following Künneth formula:
	\begin{theorem}\label{kunneth}
		Let $X$ be a locally compact, Hausdorff, second countable space, then
		$$K_i(C^*_L(X,\A(M)))\otimes\mathbb{Q}\cong \bigoplus_{j\in\mathbb{Z}_2}K_j(X)\otimes K_{i-j}(\A(M))\otimes\mathbb{Q}.$$
	\end{theorem}
	The proof of the Künneth formula is analogous to the proof of \cite[Theorem~A.10]{guo2024hilbert}. One observes that there is a canonical inclusion $\IC_L[X]\otimes\A(M)\to\IC_L[X,\A(M)]$, and both functors
	$$X\mapsto K_*(X)\otimes K_*(\A(M))\otimes\mathbb{Q}\quad\text{and}\quad X\mapsto K_*(C^*_L(X,\A(M)))\otimes\mathbb{Q}$$
	satisfy the Mayer-Vietoris theorem. When $X$ is a point, the theorem holds since 
	$$K_*(C^*_L(X,\A(M)))\cong K_*(\A(M)).$$ 
	By using a cutting and pasting argument and the Five Lemma, the theorem holds for the case when $X$ is a compact CW-complex. Since any compact Hausdorff space can be written as an inverse limit of finite CW-complexes (the nerve complex associated with a suitable open cover), the theorem holds for all compact spaces. For a locally compact space $X$, consider its one-point compactification $X^+$. The theorem follows by the exact sequence $0\to C_0(X)\to C(X^+)\to \IC\to 0$ and the Five Lemma.
	
	\begin{proof}[Proof of Theorem \ref{BLinjective}.]
		We have the following commutative diagram:
		\[\begin{tikzcd}
			K_{i+1}(C^*_L(X))\otimes\mathbb{Q} \arrow[r, "(\beta_L)_*"] \arrow[d, "\cong"'] & K_i(C^*_L(X ,\A(M)))\otimes\mathbb{Q}\arrow[d, "\cong"] \\
			\bigoplus_{j\in\mathbb{Z}_2}K_{i-j+1}(C^*_L(X))\otimes K_j(\S)\otimes\mathbb{Q}  \arrow[r, "1\otimes(\beta_{x_0})_*"] & \bigoplus_{j\in\mathbb{Z}_2}K_{i-j}(C^*_L(X))\otimes K_j(\A(M))\otimes\mathbb{Q},                
		\end{tikzcd}\]
		where the vertical maps are given by the Künneth formula, and the bottom horizontal map is given by the Bott map $\beta_{x_0}:\S\to\A(M)$ associated with the base point $x_0\in M$.
		Since $M$ is admissible, we conclude that the bottom map is injective by Proposition~\ref{bottinjective}. This implies that $(\beta_L)_*$ is a $\mathbb{Q}$-injection, which completes the proof.
	\end{proof}
	
	\subsection{The twisted assembly map}
	The evaluation homomorphism
	$$ev :C^*_{L}(P_d(X),\A(M))\to C^*(P_d(X),\A(M))$$
	defined by $ev(g)=g(0)$, induces a homomorphism on $K$-theory:
	$$ev_*:K_*(C^*_{L}(P_d(X),\A(M)))\to K_*(C^*(P_d(X),\A(M))).$$
	In this section, our goal is to prove that it is an isomorphism.
	\begin{theorem}\label{twisted Baum-Connes}
			Let \(X\) be a proper metric space with bounded geometry that admits a coarse embedding into an admissible \(\ell^p\)-Busemann space \(\mathcal M\). Then the $K$-theoretic evaluation homomorphism
		$$ev_*:K_*(C^*_{L}(P_d(X),\A(M)))\to K_*(C^*(P_d(X),\A(M)))$$
		is an isomorphism.
	\end{theorem}
	\begin{definition}
		Let $O$ be an open subset of $M\times\mathbb R_+$. Define $\IC[P_d(X),\mathcal A(M)]_O$ to be the subalgebra of $\IC[P_d(X),\mathcal A(M)]$ generated by all $T\in\mathbb C[P_d(X),\mathcal A(M)]$ such that
		$$\textup{supp}(T(x,y))\subseteq O$$
		for all $x,y\in N_d$.
		
		Define $\IC_L[P_d(X),\mathcal A(M)]_O$ to be the subalgebra of $\IC_L[P_d(X),\mathcal A(M)]$ generated by all $g\in\mathbb C_L[P_d(X),\mathcal A(M)]$ such that
		$$\textup{supp}((g(t))(x,y))\subseteq O$$
		for all $x,y\in N_d$ and $t\in\mathbb R_+$.
		Define $C^*(P_d(X),\mathcal A(M))_O$ and $C^*_L(P_d(X),\mathcal A(M))_O$ to be the norm closures of $\IC[P_d(X),\mathcal A(M)]_O$ and $\IC_L[P_d(X),\mathcal A(M)]_O$, respectively.
	\end{definition}
	
	\begin{lemma}\label{slice}
		Let $O$ be an open subset of $M\times\mathbb R_+$. Suppose $O=\bigsqcup_{i\in I}O_i$ is a disjoint union of $\{O_i\}$ in $M\times\mathbb R_+$ such that:
		\begin{enumerate}
			\item there exists $s>0$, such that $d(O_i,O_j)>s$ for $i\neq j$.
			\item there exists $r>0$ satisfying that for each $i$, there is an $x_i\in X$ with
			$$O_{i}\subseteq B(f(x_i),r).$$
		\end{enumerate} 
		Then $$ev_*:\lim_{d\to\infty}K_*\left(C^*_L(P_d(X),\mathcal A(M))_O\right)\to\lim_{d\to\infty}K_*\left(C^*(P_d(X),\mathcal A(M))_O\right)$$
		is an isomorphism.	
	\end{lemma}
	
	Let $X_R^{i}$ be the set $\{x\mid x\in X, d(x,x_i)\leq R\}$. Define
	$$\mathbb{A}(X_R)=\left\{(T_i)_{i\in I}\in\prod_{i\in I}\mathbb{C}_f(X_R^{i})\otimes\mathcal A(M)_{O_{i}} \,\Big|\, \sup_i \textup{Prop}(T_i)< \infty\right\}$$
	
	Similarly, define $\mathbb{A}_L(X_R)$ to be the $C^*$-subalgebra of $\prod_{i\in I}\mathbb{C}_L(X_R^{i})\otimes\mathcal A(M)_{O_{i}}$ generated by all elements $(b_i)_{i\in I}$ such that
	\begin{itemize}
		\item[(1)] the map $t\mapsto\prod b_i(t)$ is uniformly norm-continuous,
		\item[(2)] there exists a bounded function $c:\mathbb R_+\to\mathbb R_+$ with $\lim_{t\to\infty}c(t)=0$ such that $(b_i(t))(x,y)=0$ whenever $d(x,y)>c(t)$ for all $i\in I$, $x,y\in X_i$, and $t\in\mathbb R_+$.
	\end{itemize}
	Define $A^*(X_R)$ and $A_L^*(X_R)$ to be the norm closures of $\mathbb{A}(X_R)$ and $\mathbb{A}_L(X_R)$, respectively.
	\begin{proposition}\label{sliceiso}
		$$C^*(P_d(X),\mathcal A(M))_O\cong\lim\limits_{R\to\infty}A^*(X_R)$$
		$$C_{L}^*(P_d(X),\mathcal A(M))_O\cong\lim\limits_{R\to\infty}A_L^*(X_R)$$
	\end{proposition}
	\begin{proof}
		Since $O=\bigsqcup_{i\in I}O_{i}$, for each $T\in \mathbb{C}[P_d(X),\mathcal A(M)]_{O}$, there is a corresponding decomposition of the operator according to the support:
		$$T=\bigoplus_{i\in I}T_{i},$$
		where
		$$T_{i}(x,y)=T(x,y)|_{O_{i}}\in\mathcal A(M)_{O_{i}}.$$
		Next we prove $\bigoplus_{i\in I}T_{i}\in A^*(X_R)$ for some $R>0$. Let $T_{i}(x,y)\neq 0$, $\operatorname{supp}(T(x,y))
		\subseteq B(f(x),R_0)$ for some $R_0>0$. Since $O_{i}\subseteq B(f(x_i),r)$, we have 
		$$d_{M}(f(x),f(x_i))\leq r+R_0.$$
		Since $f$ is a coarse embedding, there is an $r'>0$ such that $d_{X}(x,x_i)\leq r'$. Since $T$ has finite propagation, there is an $r'' > 0$ such that $d_{X}(x,y)\leq r''$, so $d_{X}(x_i,y)\leq r''+r'$. Thus $y\in X^i_{r''+r'}$. Let $R=r''+r'$, which does not depend on $i\in I$. We have proved $\bigoplus_{i\in I}T_{i}\in \mathbb{A}(X_R)$. 
		One can check that the correspondence $T\mapsto\bigoplus_{i\in I}T_{i}$ extends to a $*$-isomorphism between the two $C^*$-algebras. 
	\end{proof}
	
	\begin{proposition}The evaluation map induces an isomorphism
		$$K_{*}(A_L^*(X_R)) \xrightarrow{ev_*}K_{*}(A^*(X_R))$$
	\end{proposition}
	\begin{proof}
		Since $(X_R^i)_{i\in I}$ is strongly Lipschitz homotopy equivalent to $(\{x_i\}_{i\in I})$, the proof is similar to \cite[Lemma~6.4]{yu2000coarse} by an Eilenberg swindle argument.
	\end{proof}
	The following lemma allows a Mayer-Vietoris argument, the proof is essentially the same as \cite[Lemma~6.11]{guo2024hilbert} or \cite[Lemma~6.12]{guo2024bott}, so we omit it.
	\begin{lemma}(c.f. \cite[Lemma 6.11]{guo2024hilbert})\label{glueopenset}
		Let $O^{(1)}$ and $O^{(2)}$ be open sets of $M\times\mathbb R_+$. Then
		\begin{itemize}
			\item[(1)] $C^*(P_d(X),\mathcal A(M))_{O^{(1)}}+C^*(P_d(X),\mathcal A(M))_{O^{(2)}}=C^*(P_d(X),\mathcal A(M))_{O^{(1)}\cup O^{(2)}}$
			\item[(2)] $C^*(P_d(X),\mathcal A(M))_{O^{(1)}}\cap C^*(P_d(X),\mathcal A(M))_{O^{(2)}}=C^*(P_d(X),\mathcal A(M))_{O^{(1)}\cap O^{(2)}}$
			\item[(3)] $C^*_L(P_d(X),\mathcal A(M))_{O^{(1)}}+C^*_L(P_d(X),\mathcal A(M))_{O^{(2)}}=C^*_L(P_d(X),\mathcal A(M))_{O^{(1)}\cup O^{(2)}}$
			\item[(4)] $C^*_L(P_d(X),\mathcal A(M))_{O^{(1)}}\cap C^*_L(P_d(X),\mathcal A(M))_{O^{(2)}}=C^*_L(P_d(X),\mathcal A(M))_{O^{(1)}\cap O^{(2)}}$
		\end{itemize}
	\end{lemma}
	Finally, we can prove the main theorem of this section.
	\begin{proof}[Proof of Theorem \ref{twisted Baum-Connes}]
		Since $X$ has bounded geometry, for any $s>0$ there exists $N>0$ and a decomposition
		$$X=\bigsqcup_{n=1}^{N}X_n$$
		such that $d(x,y)>s$ for each $x,y\in X_n$ with $x\neq y$. As a consequence, for each $R>0$, the family $\{B(f(x),R)\}_{x\in X}$ also has a decomposition
		$$\bigcup_{x\in X}B(f(x),R)=\bigcup_{n=1}^{N}\bigsqcup_{x\in X_n}B(f(x),R).$$ Define
		$$O_n(R)=\bigsqcup_{x\in X_n}B(f(x),R).$$
		By Lemma \ref{slice}, the evaluation map induces an isomorphism on $K$-theory:
		$$ev_*:\lim_{d\to\infty}K_*\left(C^*_L(P_d(X),\mathcal A(M))_{O_n(R)}\right)\xrightarrow{\cong}\lim_{d\to\infty}K_*\left(C^*(P_d(X),\mathcal A(M))_{O_n(R)}\right).$$
		Define
		$$B(R)=\bigcup_{x\in X}B(f(x),R)=\bigcup_{n=1}^{N}O_n(R).$$
		Using Lemma \ref{glueopenset} and a Mayer-Vietoris argument $N$ times, we have that
		$$ev_*:\lim_{d\to\infty}K_*\left(C^*_L(P_d(X),\mathcal A(M))_{B(R)}\right)\xrightarrow{\cong}\lim_{d\to\infty}K_*\left(C^*(P_d(X),\mathcal A(M))_{B(R)}\right)$$
		is an isomorphism for each $R>0$. Consider the following commuting diagram:
		$$\begin{tikzcd}
			\lim\limits_{d\to\infty}K_*\left(C^*_L(P_d(X),\mathcal A(M))\right) \arrow[r,"ev_*"] \arrow[d,"\cong"'] & 
			\lim\limits_{d\to\infty}K_*\left(C^*(P_d(X),\mathcal A(M))\right) \arrow[d,"\cong"] \\
			\lim\limits_{d\to\infty}\lim\limits_{R\to\infty}K_*\left(C^*_L(P_d(X),\mathcal A(M))_{B(R)}\right) \arrow[r,"ev_*"] \arrow[d,"\cong"'] & 
			\lim\limits_{d\to\infty}\lim\limits_{R\to\infty}K_*\left(C^*(P_d(X),\mathcal A(M))_{B(R)}\right) \arrow[d,"\cong"] \\
			\lim\limits_{R\to\infty}\lim\limits_{d\to\infty}K_*\left(C^*_L(P_d(X),\mathcal A(M))_{B(R)}\right) \arrow[r,"ev_*","\cong"']           & 
			\lim\limits_{R\to\infty}\lim\limits_{d\to\infty}K_*\left(C^*(P_d(X),\mathcal A(M))_{B(R)}\right)         
		\end{tikzcd}$$
		The vertical maps are isomorphisms and the bottom horizontal map is an isomorphism as shown above. This completes the proof of Theorem \ref{twisted Baum-Connes}.
	\end{proof}
	
	Now we prove the main result, Theorem \ref{mainthm}.
	\begin{proof}[Proof of Theorem \ref{mainthm}.]
		We have the following commutative diagram:
		\[\begin{tikzcd}
			K_{*+1}(C^*_L(P_d(X))) \arrow[r, "ev_*"] \arrow[d, "(\beta_L)_*"] & K_{*+1}(C^*(P_d(X))) \arrow[d, "\beta_*"]   \\
			K_*(C^*_L(P_d(X),\mathcal A(M)))\arrow[r, "ev_*"]& K_*(C^*(P_d(X),\mathcal A(M)))
		\end{tikzcd}\]
		We see that $(\beta_L)_*$ is $\mathbb{Q}$-injective by Proposition \ref{BLinjective}, and the bottom map is an isomorphism by Theorem \ref{twisted Baum-Connes}. This implies the assembly map is $\mathbb{Q}$-injective.
	\end{proof}
	
	\bibliographystyle{abbrv}
	\bibliography{lpBusemann}
	
	\noindent  Chao Tong  \\
	Chongqing University (at Huxi Campus),\\
	Chongqing 401331, P. R. China.\\
	E-mail: \url{20131989@cqu.edu.cn}\\
	
	\noindent  Xianjin Wang \\
	Chongqing University (at Huxi Campus),\\
	Chongqing 401331, P. R. China.\\
	E-mail: \url{xianjinwang@cqu.edu.cn}\\
	
\end{document}